\documentclass[12pt]{article}
\usepackage{amssymb,amsmath,amsthm,latexsym}
\usepackage{dsfont}
\usepackage{amssymb}
\usepackage{indentfirst}
\usepackage{hyperref}
\begin{document}
\title{A new proof of the Larman--Rogers upper bound for the chromatic number of the Euclidean space}
\author{R.I. Prosanov\thanks{The author is supported in part by the grant from Russian Federation President Program of Support for Leading Scientific Schools 6760.2018.1.}}
\date{}
\AtEndDocument{\bigskip{\footnotesize%
  \textsc{Universit\'{e} de Fribourg, Chemin du Mus\'{e}e 23, CH-1700 Fribourg, Switzerland} \par
  \textsc{Moscow Institute Of Physics And Technology, Institutskiy per. 9, 141700, Dolgoprudny, Russia} \par
  \textit{E-mail}: \texttt{rprosanov@mail.ru}
}}
\date{}
\maketitle

\renewcommand{\refname}{References}
\renewcommand{\proofname}{Proof}
\renewcommand{\le}{\leqslant}
\renewcommand{\leq}{\leqslant}
\renewcommand{\ge}{\geqslant}
\renewcommand{\geq}{\geqslant}
\renewcommand{\mathds}{\mathbb}
\newcommand{\R}{\mathbb{R}}

\begin{abstract}

The \textit{chromatic number} $\chi(\R^n)$ of the Euclidean space $\R^n$ is the smallest number of colors sufficient for coloring all points of the space in such a way that any two points at the distance 1 have different colors. In 1972 Larman--Rogers proved that $\chi(\mathds{R}^n) \leq (3 + o(1))^n$. We give a new proof of this bound.

\end{abstract}

\section{Introduction}

The \textit{chromatic number} $\chi(\R^n)$ of the Euclidean space $\R^n$ is the smallest number of colors sufficient for coloring all points of this space in such a way that any two points at the distance 1 have different colors. This problem was initially posed by Nelson for $n=2$ (see the history of this problem in \cite{Ra2}, \cite{Rai2}, \cite{Rai1}, \cite{So}).

The exact value of $\chi(\R^n)$ is not known even in the planar case. The best known bounds are $$5 \leq \chi(\mathds{R}^2)
\leq 7.$$ See~\cite{So} for the upper bound and~\cite{dG} for the lower one. In the case of growing $n$ we have 

\begin{equation}
(1.239 + o(1))^n \leq \chi(\mathds{R}^n) \leq (3 + o(1))^n.
\end{equation}

The lower bound is due to Raigorodskii \cite{Ra1} and the upper bound is due to Larman and Rogers \cite{Lar}.

The proof of Larman and Rogers is based on a hard theorem due to Butler~\cite{Bu} and on a result of Erd\H os and Rogers about coverings of $\R^n$ with translates of a convex body. In this paper we present a new proof that does not use neither of them. Instead we adapt the approach developed by Marton Nasz\'{o}di in~\cite{Na}. It connects geometrical covering problems with coverings of finite hypergraphs. An advantage of this approach in contrast to the previous one is that it could be turned into an algorithmic one. Indeed, the original paper of Erd\H os and Rogers used probabilistic arguments. Therefore, the same is true for the Larman and Rogers proof. We will rely only on a theorem by Johnson, Lov\'asz and Stein that establishes a connection between the fractional covering number of a hypergraph and its integral covering number. The proofs of this theorem are quite easy and provide an algorithm which constructs an economical covering based on an optimal fractional covering. A problem of finding an optimal fractional covering is a problem of linear programming.

The author developed this method further in~\cite{Pr}, where new upper bounds for chromatic numbers of spheres were obtained.

We obtain the upper bound in (1) from a slightly more general result, stated in the next section. It is motivated by the following generalization of $\chi(\R^n)$. Let $K$ be a convex centrally-symmetric body. Consider the space $\R^n$  with the norm determined by $K$. Let $\chi(\mathds{R}_K^n)$ be the chromatic number of this normed space. If $B^n$ is the usual unit ball in $\R^n$, then we have $\chi(\R^n) = \chi(\mathds{R}_{B^n}^n)$. In 2008 Kang and F\"{u}redi \cite{Fu} obtained an upper bound for an arbitrary $K$: $$\chi(\mathds{R}_K^n) \leq (5 + o(1))^n.$$ In 2010 Kupavskii \cite{Ku} improved it to $$\chi(\mathds{R}_K^n) \leq (4 + o(1))^n.$$ No exponential lower bounds are known for the general case, although such bounds are known to hold for the case of $l_p$-norms \cite{Ra_lq}.

In the statement of our main theorem we give an upper bound for $\chi(\mathds{R}_K^n)$ in terms of another quantity: the \textit{tiling parameter} of $K$ (see below for the precise definition). It is of interest to investigate this quantity on its own. In particular, the \textit{lattice tiling parameter} measures the Banach-Masur distance from a centrally-symmetric convex body to a closest parallelohedron. We show that any progress on bounding the tiling parameter will lead to a progress on bounding $\chi(\mathds{R}_K^n)$. 

From the paper of Butler a bound on the lattice tiling parameter of $K=B^n$ can be established. In Section 3 we demonstrate that for our generalization a similar bound can be obtained without any efforts and then applied to chromatic numbers.

This paper is organized as follows. In Section 2 we give all necessary definitions and formulate the main result of this paper. In Section 3 we deduce the upper bound in (1) from this result.

\section{The main result}

\subsection{Multilattices and tiling parameters}

Let $\Omega$ be a lattice. A \textit{multilattice} is a union $\Phi = \cup_{i=1}^{q} (\Omega+x_i)$ of translates of $\Omega$ by a finite number of vectors. A lattice can be considered as a multilattice with $q=1$. The number $q$ in the definition of $\Phi$ is denoted by $q(\Phi)$. A tiling $\Psi$ of the space $\R^n$ by convex polytopes is called \textit{associated} with the multilattice $\Phi$ if there is a bijection between polytopes of $\Psi$ and points in $\Phi$ such that every point $x$ is contained in the interior of the corresponding polytope $\psi_x$.

Let $K$ be a bounded closed centrally-symmetric convex body. The \textit{tiling parameter} is $$\gamma(K, \Phi, \Psi)=\inf \{\beta/\alpha | \forall x \in \Phi, \alpha K + x \subset \psi_x \subset \beta K + x \}.$$ Define $$\gamma(K, k)=\inf_{\Phi, \Psi :~ q(\Phi) \leq k} \gamma(K, \Phi, \Psi),$$ where the infimum is taken over all multilattices $\Phi$ with $q(\Phi) \leq k$ and all tilings associated with them. We call $\gamma(K,1)$ to be the \textit{lattice tiling parameter} of $K$.

Our main result is

\vskip+0.2cm
{\bfseries Theorem 1.} \textit {We have $$\chi(\mathds{R}_{K}^n) \leq (1 + \gamma(K, k))^n\bigg[n\ln n+n\ln \ln n+2\ln k+2n\Big(1+\ln \big(2 \gamma(K, k)\big)\Big)\bigg].$$ In particular, if $K_n$ is a sequence of bodies and $k_n$ is a sequence of positive numbers such that for some absolute constant $c$ we have $k_n\leq n^{cn}$, then  $$\chi(\mathds{R}_{K_n}^n) \leq (1 + \gamma(K_n, k_n)+o(1))^n.$$
}
\vskip+0.2cm

\subsection{Preliminaries with fractional coverings}

Let $Z$ be a set, $\mathcal{F}$ be a family of its subsets, and $Y \subseteq Z$. By the \emph{covering number} $\tau(Y,\mathcal{F})$ denote the minimal cardinality of a family $\mathcal{H} \subseteq \mathcal{F}$ such that $Y$ is covered by the union of all sets $F \in \mathcal{H}$.

If $Z$ is finite, then the pair $(Z, \mathcal F)$ is a finite hypergraph. In this case a \textit{fractional covering of $Y$ by $\mathcal{F}$} is a function $\nu: \mathcal{F} \rightarrow  [0; +\infty)$ such that for all $y \in Y$ we have$$\sum\limits_{F \in \mathcal{F}: y \in F} \nu(F) \geq 1.$$ Define the \emph{fractional covering number} of $Y$:
$$\tau^*(Y,\mathcal{F})= \inf \{ \nu(\mathcal{F}): \nu{\rm~is~a~fractional~covering~of~} Y{\rm~by~}\mathcal{F} \}.$$

The following theorem establishes a connection between $\tau$ and $\tau^*$

\vskip+0.2cm
{\bfseries Theorem 2 (\cite{Jo}, \cite{Lo}, \cite{St}.)} \textit {Suppose $Z$ is a finite set and $\mathcal{F} \subseteq 2^{Z}$, then
$$
\tau(Z,\mathcal{F}) < \left(1+ \ln\left(\underset{F \in \mathcal{F}}{\max} (|F|)\right)\right)\tau^*(Z,\mathcal{F}).$$
}
\vskip+0.2cm

\subsection{Proof of Theorem 1}

In what follows, all distances are calculated with respect to the norm, determined by~$K$. For $0<\mu<1$ define $$\mu\Psi=\bigcup\limits_{\psi_x \in \Psi} \left(\mu (\psi_x - x) + x\right).$$


Fix $\varepsilon > 0$. Choose a pair $(\Phi, \Psi)$ such that $$\gamma(K, \Phi, \Psi) < \gamma(K, k) + \varepsilon.$$ Let $\alpha, \beta$ be numbers such that $$\gamma=\beta/\alpha < \gamma(K, \Phi, \Psi) + \varepsilon$$ and for all $x  \in \Phi$ we have $$\alpha K + x \subset {\rm int}(\psi_x),$$ $$\psi_x \subset {\rm int}(\beta K + x).$$ Since $\psi_x$ is contained in ${\rm int}(\beta K+x)$, we see that the diameter of $\psi_x$ is strictly less than~$2\beta$. Let $$\mu=\alpha/(\alpha+\beta).$$ Then for all $x$ the polytope $\psi_x$ does not contain a pair of points at distance $2\beta\mu$. 

We show that for all $x, y \in \Phi$, $x \neq y$, the distance between $\mu\psi_x$ and $\mu\psi_y$ is greater than $2\beta\mu$. It is sufficient to consider only such polytopes that share a common face in some dimension. Since $\psi_x$ and $\psi_y$ are convex and share some $k$-dimensional face, there is a hyperplane containing this face and separating $\psi_x$ and $\psi_y$. Let $l_x$ and $l_y$ be the distances from $x$ and $y$ to this hyperplane. The distance between $\mu \psi_x$ and $\mu \psi_y$ is greater than the distance between the images of this hyperplane under homothety with center $x$ and homothety with center $y$. Since $\alpha K + x \subset {\rm int}(\psi_x)$, this distance is $$(l_x+l_y)(1-\mu) > 2\alpha(1-\mu)=2\beta\mu.$$ Therefore, the set $\mu\Psi$ does not contain a pair of points at the distance $2\beta\mu$ and we can color it with one color.

Next, we cover $\R^n$ by the copies of $\mu \Psi$. This set is a disjoint union of several convex bodies. Hence, typical covering results (like in~\cite{ER}) can not be applied to it. Now we show how to overcome this difficulty.

Let $\Omega$ be the base lattice of $\Phi$. Consider the torus $T^n=\mathds{R}^n/\Omega$. Let $\tilde x_i$ be the projections onto $T^n$ of the translation vectors $x_i$ of the lattice $\Omega$ in the multilattice $\Phi$ and $\tilde X$ be their union. The tiling $\Psi$ is periodical over the lattice $\Omega$, hence we can define its projection $\tilde{\Psi}$, which is a tiling of $T^n$ associated to the set $\tilde X$. 

We will cover $T^n$ by less than $$(1 + \gamma)^n(n\ln n+n\ln \ln n+2\ln k+2n(1+\ln 2 \gamma))$$ translates of $\mu \Tilde {\Psi}$.

We need the following lemma.

\vskip+0.2cm
{\bfseries Lemma 1.} {\it Fix $0 < \delta < 1$. Let $\mathcal{F}$ and $\mathcal{F'}$ be the families of translates of the sets $\mu \Tilde {\Psi}$ and $\mu (1- \delta) \Tilde {\Psi}$ by all points of $T^n$. Suppose $\Lambda \subset T^n$ is a finite point set of maximal cardinality such that $\frac{\alpha \mu \delta }{2} K + \Lambda$ is a packing of the bodies $\frac{\alpha \mu \delta }{2} K$. Then $\tau(T^n,\mathcal{F}) \leq \tau(\Lambda,\mathcal{F'})$.}
\vskip+0.2cm

\begin{proof}
Since the cardinality of $\Lambda$ is maximal, then $\alpha \mu \delta K + \Lambda$ is a covering of $T^n$.

Let $Y = \{ y_j,~j = 1, \dots, m \} \subset T^n$ be a point set such that $\mu (1- \delta) \Tilde {\Psi} + Y$ covers $\Lambda$. We show that  $\mu \Tilde {\Psi}+Y$ covers $T^n$.

Let $t \in T^n$ be an arbitrary point. Since $\alpha \mu \delta K + \Lambda$ is a covering of $T^n$, then there exists $\lambda \in \Lambda$ such that $\alpha \mu \delta K+\lambda$ contains $t$. There also exists $j$ such that $$\lambda \in ((\mu (1- \delta) \tilde {\psi}_i+y_j))$$ for some $i$. Since for all $i$ we have $\alpha K \subset {\rm int}(\tilde {\psi}_i - \tilde {x}_i)$, we obtain $$t \in \mu \delta \tilde {\psi}_i - \tilde {x}_i + \lambda \subset ( (\mu \delta \tilde {\psi}_i - \tilde {x}_i) + (\mu (1-\delta) \tilde {\psi}_i - \tilde {x}_i) ) + \tilde {x}_i + y_j \subset \mu \tilde {\psi}_i - \tilde {x}_i + \tilde {x}_i +y_j = \mu \tilde {\psi}_i+y_j.$$ The proof is complete.

\end{proof}

Consider $\mathcal{F}$, $\mathcal{F'}$ and $\Lambda$ as in the notation of Lemma 1. Define $$\mathcal E=\{ \Lambda \cap F: F \in \mathcal F' \}.$$

Then $(\Lambda, \mathcal E)$ is a finite hypergraph and $\tau(\Lambda,\mathcal{F'})=\tau(\Lambda,\mathcal{E}).$

From Lemma 1 and Theorem 2 it follows that $$\tau(T^n,\mathcal{F}) \leq \tau(\Lambda,\mathcal{E}) \leq \left(1+ \ln\left(\underset{E \in \mathcal{E}}{\max} (|E|)\right)\right)\tau^*(\Lambda,\mathcal{E}).$$

We want to bound $\tau^*(\Lambda,\mathcal{E})$. By $\sigma$ denote the usual measure on $T^n$ induced by the Lebesgue measure on $\R^n$ and scaled in such a way that $\sigma(T^n)=1$. 

For $\lambda \in \Lambda$ and $E \in \mathcal E$ define $$S(\lambda) =\{t \in T^n: \lambda \in \mu (1- \delta) \Tilde {\Psi} + t \}, $$ $$S(E)=\{t \in T^n: E \subseteq \mu (1- \delta) \Tilde {\Psi} + t\}.$$

The sets $S(\lambda)$, $S(E)$ are measurable. Moreover, for every $\lambda \in \Lambda$, $$\sigma(S(\lambda)) =  \sigma(\mu (1- \delta) \Tilde {\Psi}),$$ $$\sigma(S(\lambda))=\sum\limits_{\lambda \in E}\sigma(S(E)),$$ $$\sum\limits_{E\in \mathcal E}\sigma(S(E)) = \sigma(T^n)=1.$$

Define $\nu: \mathcal E \rightarrow [0; +\infty)$ as $$\nu(E)=\frac{\sigma(S(E))}{\sigma(\mu (1- \delta) \Tilde {\Psi})}.$$

Then it is a fractional covering of $\Lambda$ by $\mathcal E$ and $$\tau^*(\Lambda, \mathcal E) \leq \sum\limits_{E \in \mathcal E} \nu(E) = \frac{1}{\sigma(\mu (1- \delta) \Tilde {\Psi})}.$$


Now we bound $\max\limits_{E \in \mathcal E}|E|=\max\limits_{F' \in \mathcal {F'}}| \Lambda \cap F' |$. 

Recall that $\frac{\alpha \mu \delta }{2} K + \Lambda$ is a packing of bodies $\frac{\alpha \mu \delta }{2} K$. If $\lambda \in F'=\mu (1- \delta) \Tilde {\Psi}+t,$ then $\lambda \in \mu (1- \delta) \tilde {\psi}_i+t$ for some $i$. 

Since $\alpha K \subset {\rm int}(\tilde {\psi}_i - \tilde {x}_i)$, we have $$\frac{\alpha \mu \delta }{2} K + \lambda \subseteq \left(\mu \frac\delta2 \tilde {\psi}_i - \tilde {x}_i\right) +
(\mu (1-\delta) \tilde {\psi}_i - \tilde {x}_i) + \tilde {x}_i + t \subseteq \mu \left(1-\frac\delta2\right) \tilde {\psi}_i+t \subseteq \mu \tilde {\psi}_i+t.$$ Now we can compare the volumes and get the bound for $| \Lambda \cap F' |$. Since every $\tilde {\psi}_i$ is contained in $\beta K+\tilde {x}_i$, we get ${\rm vol}(\tilde {\psi}_i) \leq \beta^n {\rm vol}(K)$ and $$| \Lambda \cap F' | \leq \frac { {\rm vol}(\mu \Tilde {\psi}) } { {\rm vol}(\frac{\alpha \mu \delta }{2} K) } \leq \sum_{i=1}^k \frac { {\rm vol}(\mu \phi_i) } { {\rm vol}(\frac{\alpha \mu \delta }{2} K) } \leq \frac{ k\mu^n\beta^n {\rm vol}(K) } {\mu^n\alpha^n(\delta/2)^n {\rm vol}(K) }=k(2\gamma/\delta)^n.$$

Finally, we obtain $$\tau(T^n,\mathcal{F}) \leq \left(1+ \ln\left(\underset{F' \in \mathcal{F'}}{\max} (| \Lambda \cap F' |)\right)\right)\tau^*(T^n,\mathcal{F'}) \leq \frac{(1+n \ln(2\gamma/\delta)+\ln k)(1+\gamma)^n}{(1-\delta)^n}.$$

Now we take $\delta = \frac{1}{2n \ln n}$ and use (for arbitrary large $n$) $$
\Big(1 - \frac{1}{2n \ln n}\Big)^{-n} \leq \exp \Big(\frac{1}{\ln n}\Big) \leq 1+\frac{2}{\ln n}.
$$

Thus we have $$\tau(T^n,\mathcal{F}) \leq (1 + \gamma)^n\left(1+\frac{2}{ \ln n}\right)\left(1+
n \ln{((4n) (\ln n) (\gamma))} + \ln k\right) \leq
$$
$$
\leq (1 + \gamma)^n\left(1 + n\ln n+n\ln \ln n
+ 2 n\left(1 + \ln 2 + \frac {\ln 4}{\ln n} + \frac {\ln {\ln n}}{\ln n}\right) +
\frac{2}{\ln n} + \right.
$$
$$
+ \left. \left(1+ \frac{2}{\ln n}\right)(n \ln \gamma + \ln k) \right) \leq
$$
$$
\leq (1 + \gamma)^n(n\ln n+n\ln \ln n+2\ln k+2n + 2n \ln (2 \gamma))\leq$$ $$\leq (1
+ \gamma(K, k) + \varepsilon)^n\left[n\ln n+n\ln \ln n+2\ln k+2n(1+\ln (2
(\gamma(K, k)+\varepsilon)))\right].$$

This inequality holds for every $ \varepsilon> 0 $. This completes the proof of Theorem 1.

\section{Chromatic number for the Euclidean metric}

In the paper \cite{Lar}, Larman and Rogers proved that for a Euclidean ball $B^n$, then the lattice tiling parameter $\gamma(B^n, 1) \leq 2+o(1)$ as $n \rightarrow \infty$. They used a theorem due to Butler \cite{Bu}. We need some notation to state the Butler result.

Let ${\mathcal K}=K+\Omega$ be a system of translates of $K$ by the vectors of the lattice $\Omega$, $\xi_1=\xi_1({\mathcal K})$ be the infimum of the positive numbers $\xi$ such that the system $\xi{\mathcal K}$ is a covering of $\R^n$, and $\xi_2=\xi_2({\mathcal K})$ be the supremum of the positive numbers $\xi$ such that $\xi{\mathcal K}$ is a packing in $\R^n$. 

Denote $\xi({\mathcal K})=
\xi_1({\mathcal K})/\xi_2({\mathcal K})$. Consider $\tilde{\gamma}(K)=\inf\limits_{\mathcal K}\xi({\mathcal K})$, where the infimum is over the set of all lattices in $\R^n$. By $DK$ denote the difference body of $K$, i.e. ${DK = \{x-y |~x, y \in K \}}$. 

\vskip+0.2cm	
{\bfseries Theorem 3 (Butler, \cite{Bu}).} \textit {Let $K$ be a bounded convex body in $\R^n$, then there exists an absolute constant $c$ such that
$$
\tilde{\gamma}(K) \leq \left[\frac{{\rm {\rm vol}}(DK)}{{\rm {\rm vol}}(K)}n^{\log_2(\ln n)+c}\right]^{1/n}.
$$}
\vskip+0.2cm

If $K$ is centrally symmetric, then we get $\tilde{\gamma}(K) \leq 2 + o(1)$. It is easy to see that if $K=B^n$, then $\gamma(B^n, 1) \leq \tilde{\gamma}(B^n)$. 

Indeed, let $\Omega$ be a lattice such that $\xi({\mathcal K}) < \tilde{\gamma}(B^n) + \varepsilon$ and $\Psi$ be a Voronoi tiling which corresponds to $\Omega$. Then for all $x \in \Omega$ we get $$K + x \subset \psi_x \subset \xi({\mathcal K}) K + x.$$ Since $\varepsilon$ is arbitrary close to zero, we obtain our inequality. Unfortunately, if $K$ is not a Euclidean ball, then Voronoi polytopes might be nonconvex and the locus of the points that have equal distances to a pair of given points might have nonzero measure. Therefore, the problem of bounding $\gamma(K, 1)$ becomes much harder.

The proof of Theorem 3 is quite nontrivial. But the problem of bounding of our generalized tiling parameter instead of the lattice one is much easier. First, we show that for some $k$, $\gamma(B^n, k) \le 2$. 

Let $\Omega$ be a lattice such that ${\mathcal K}=B^n+\Omega$ is a packing. We claim that there is some multilattice $\Phi$ with the base lattice $\Omega$ such that ${\mathcal B^n}=K+\Omega$ is a packing and $2{\mathcal K}=2B^n+\Omega$ is a covering. By $T^n$ denote the torus $\mathds{R}^n/\Omega$. Choose a set $Y=\{y_i\}$ of the maximal cardinality such that $$(B^n+y_i)\cap (B^n+y_j)={\varnothing}, ~~~ \forall i\neq j.$$ For all $x \in T^n$ there exists $i$ such that $\|x-y_i\|_K<2$. Otherwise, $(B^n+x)\cap(B^n+y_i)={\varnothing}$ which implies that $Y$ does not have the maximal cardinality. Therefore, we have proved that $\bigcup\limits_i 2B^n+y_i$ covers $T^n$. Hence, we can take the multilattice $\Phi=\Omega+Y$.

Now associate to it the Voronoi tiling $\Psi$ of the point set $\Omega+Y$. Then in turn $$\gamma(B^n, k) \leq \gamma(B^n, \Phi, \Psi)   \leq 2.$$

Now we supply an upper bound on $k$. Let $B^n$ be inscribed into a cube $C$ with the side length 2. The edges of $C$ generate a lattice $\Omega$ in $\R^n$ such that ${\mathcal K}=B^n+\Omega$ is a packing. We can bound $k$ using volumes $$k \leq \frac{{\rm vol}(T)}{{\rm vol}(B^n)} \leq n^{cn}.$$

Finally, we can apply Theorem 1 and get the upper bound for (1).

\vskip+0.2cm
\textbf{Acknowledgements.} The author is grateful to A. M. Raigorodskii and A. B. Kupavskii for their constant attention to this work and for useful remarks.
\vskip+0.2cm

\bibliographystyle{abbrv}
\bibliography{lar-rod}

\end{document}